\documentclass[a4paper,12pt]{amsart}

\usepackage{amsmath}
\usepackage{amssymb}
\usepackage{mathrsfs}
\usepackage{ifthen}
\usepackage{graphicx}
\usepackage[T1]{fontenc} 

\def\switchlinenumbers{\@ifstar
    {\let\makeLineNumberOdd\makeLineNumberRight
     \let\makeLineNumberEven\makeLineNumberLeft}%
    {\let\makeLineNumberOdd\makeLineNumberLeft
     \let\makeLineNumberEven\makeLineNumberRight}%
    }

\def\setmakelinenumbers#1{\@ifstar
  {\let\makeLineNumberRunning#1%
   \let\makeLineNumberOdd#1%
   \let\makeLineNumberEven#1}%
  {\ifx\c@linenumber\c@runninglinenumber
      \let\makeLineNumberRunning#1%
   \else
      \let\makeLineNumberOdd#1%
      \let\makeLineNumberEven#1%
   \fi}%
  }

\nonstopmode \numberwithin{equation}{section}
\newtheorem*{theorem*}{Theorem}

\newtheorem{thm}{Theorem}[section]
\newtheorem{cor}[equation]{Corollary}
\newtheorem{lem}[equation]{Lemma}
\newtheorem{prop}[equation]{Proposition}

\theoremstyle{definition}
\newtheorem{defn}{Definition}[section]

\newtheorem{prob}[equation]{Problem}

\newtheorem{innercustomthm}{Theorem}
\newenvironment{customthm}[1]
  {\renewcommand\theinnercustomthm{#1}\innercustomthm}
  {\endinnercustomthm}

\newcounter{minutes}
\divide\time by 60
\newcounter{hours}
\multiply\time by 60
\addtocounter{minutes}{-\time}

\newcounter {own}
\def\theown {\thesection       .\arabic{own}}

\newenvironment{pf}[1][]{%
 \vskip 3mm
 \noindent
 \ifthenelse{\equal{#1}{}}%
  {{\slshape Proof. }}%
  {{\slshape #1.} }%
 }%
{\qed\bigskip}

\newcounter{alphabet}

\def\be{\begin{equation}}
\def\ee{\end{equation}}

\newcommand{\bee}{\begin{enumerate}}
\newcommand{\eee}{\end{enumerate}}

\newcommand{\blem}{\begin{lem}}
\newcommand{\elem}{\end{lem}}
\newcommand{\bthm}{\begin{thm}}
\newcommand{\ethm}{\end{thm}}
\newcommand{\bcor}{\begin{cor}}
\newcommand{\ecor}{\end{cor}}
\newcommand{\beg}{\begin{examp}}
\newcommand{\eeg}{\end{examp}}
\newcommand{\begs}{\begin{examples}}
\newcommand{\eegs}{\end{examples}}
\newcommand{\bdefe}{\begin{defin}}
\newcommand{\edefe}{\end{defin}}
\newcommand{\bprob}{\begin{prob}}
\newcommand{\eprob}{\end{prob}}
\newcommand{\bei}{\begin{itemize}}
\newcommand{\eei}{\end{itemize}}

\newcommand{\norm}[1]{\left\lVert#1\right\rVert}

\begin{document}

\title{Arithmetic Bohr radius of bounded linear operators}

\author{Vasudevarao Allu}
\address{Vasudevarao Allu,
Department of Mathematics,
School of Basic Sciences,
Indian Institute of Technology Bhubaneswar,
Bhubaneswar-752050, Odisha, India.}
\email{avrao@iitbbs.ac.in}

\author{Subhadip Pal}
\address{Subhadip Pal,
	Department of Mathematics,
	School of Basic Sciences,
	Indian Institute of Technology Bhubaneswar,
	Bhubaneswar-752050, Odisha, India.}
\email{subhadippal33@gmail.com}

\subjclass[{AMS} Subject Classification:]{Primary 32A05, 32A10, 46B07; Secondary 46G20}

\def\thefootnote{}
\footnotetext{ {\tiny File:~\jobname.tex,
printed: \number\year-\number\month-\number\day,
          \thehours.\ifnum\theminutes<10{0}\fi\theminutes }
} \makeatletter\def\thefootnote{\@arabic\c@footnote}\makeatother

\begin{abstract}
In this paper, we investigate the arithmetic Bohr radius of bounded linear operators between arbitrary complex Banach spaces. We establish the close connection between the classical Bohr radius and the arithmetic Bohr radius of bounded linear operators. Further, we study the asymptotic estimates of arithmetic Bohr radius for identity operator on infinite dimensional complex Banach spaces. Finally, we obtain the correct asymptotic behavior of Bohr radii of operators between sequence spaces.
\end{abstract}

\maketitle
\pagestyle{myheadings}
\markboth{Vasudevarao Allu and Subhadip Pal}{Arithmetic Bohr radius of bounded linear operators}

\section{Introduction}
In recent years, studying different multidimensional variants of Bohr's classical power series theorem have been an active field of research and it has been studied extensively in diverse articles, for example \cite{aizn-2000a,aizn-2000b,aizenberg-2001,boas-1997,boas-2000,defant-2006,defant-2003,defant-2004,Dineen-Timoney-1989,Djakov & Ramanujan & J. Anal & 2000,paulsen-2002}. Following Boas and Khavinson \cite{boas-1997}, the multidimensional Bohr radius $K_n$ of $n$-dimensional polydisk $\mathbb{D}^n=\{z=(z_1,\ldots,z_n)\in \mathbb{C}^n : \max_{1\leq j\leq n}|z_j|<1\}$ in $\mathbb{C}^n$ is the supremum of all positive numbers $r$ such that 
\begin{equation*}
	\sup_{z\in r\mathbb{D}^n}\sum_{\alpha\in \mathbb{N}^n_{0}}|c_{\alpha}z^{\alpha}| \leq \sup_{z\in \mathbb{D}^n}\bigg|\sum_{\alpha\in \mathbb{N}^n_{0}}c_{\alpha} z^{\alpha}\bigg|
\end{equation*} 
for all holomorphic functions $\sum_{\alpha\in \mathbb{N}^n_{0}}c_{\alpha} z^{\alpha}$ on $\mathbb{C}^n$. In 1914, H. Bohr \cite{Bohr-1914} himself studied and estimated the constant $K_{1}$, the classical Bohr radius and it was shown that $K_{1}$ to be $1/6$.
Later, the Bohr radius $K_1$ was further improved to $1/3$ by M. Riesz, I. Schur, and N. Wiener \cite{sidon-1927,tomic-1962}, independently. It is worth mentioning here that for $n\geq 2$, the exact values of $K_n$ are unknown and till today it is an open problem since 1997. However, Boas and Khavinson \cite{boas-1997} have established the following estimate for $K_n$:
\begin{equation*}
	\frac{1}{3}\sqrt{\frac{1}{n}}\leq K_n \leq 2\sqrt{\frac{\log n}{n}}.
\end{equation*}
The significant work by Boas and Khavinson \cite{boas-1997} became the source of inspiration for many subsequent papers on multidimensional Bohr radii connecting the asymptotic estimates of Bohr radius to different problems in Functional Analysis like geometry of Banach spaces \cite{Blasco-Collect-2017}, unconditional basis constants of spaces of polynomials \cite{aizenberg-2001}, operator algebra \cite{Dixon & BLMS & 1995,paulsen-2002} etc.\\

For the last decade, an extensive study has been carried out on multidimensional Bohr radius by replacing the scalar-valued holomorphic functions by Banach space valued ones. Blasco \cite{Blasco-OTAA-2010} has first initiated the concept of Bohr radius $K(\mathbb{D},X)$ for vector-valued holomorphic functions defined on unit disk $\mathbb{D}$ in complex plane into a complex Banach space $X$. Later, the author \cite{Blasco-Collect-2017} has extended the study to $p$-Bohr radius $K^p(\mathbb{D},X)$ for $1\leq p<\infty$. In 2012, Defant {\it et al.} \cite{defant-adv-math-2012} introduced a more general version of multidimensional Bohr radius of bounded linear operators between arbitrary complex Banach spaces. In \cite{Subhadip-Vasu-P5}, motivated by the work of Defant {\it et al.} \cite{defant-adv-math-2012} and Blasco \cite{Blasco-Collect-2017}, we have introduced the following notion of Bohr radius of bounded linear operators. 
\begin{defn}\label{Pal-Vasu-P5-defn-01}\cite{Subhadip-Vasu-P5}
	Let $1\leq q\leq \infty$, $1\leq p<\infty$, $\lambda\geq 1$, $n\in \mathbb{N}$, and let $T: X \rightarrow Y$ be a bounded linear operator between complex Banach spaces $X$ and $Y$ such that $\norm{T}\leq \lambda^{1/p}$. For $1\leq p<\infty$, the $\lambda_{p}$-Bohr radius of $T$, denoted by $K^p(B_{\ell^n_{q}}, T,\lambda)$, is the supremum of all $r\geq 0$ such that for all holomorphic functions $f(z)=\sum_{\alpha\in \mathbb{N}^n_{0}}x_{\alpha}z^{\alpha}$ on $B_{\ell^n_{q}}$, we have 
	\begin{equation}\label{Pal-Vasu-P5-eqn-int-04}
		\sup_{z\in rB_{\ell^n_{q}}}\sum_{\alpha\in \mathbb{N}^n_{0}}\norm{T(x_{\alpha})z^{\alpha}}^{p}_{Y}\leq \lambda\norm{\sum_{\alpha\in \mathbb{N}^n_{0}}x_{\alpha}z^{\alpha}}^p_{X}.
	\end{equation}
\end{defn}
\noindent
We write $K^p(B_{\ell^n_q},T)$ for $\lambda=1$. If $T$ is the identity map on $X$ then we use the notation $K^p(B_{\ell^n_q},X,\lambda)$ and $K^p(B_{\ell^n_q},X)$. For $X=\mathbb{C}$, we write $K^p(B_{\ell^n_q},\lambda)$ for $K^p(B_{\ell^n_q},\mathbb{C},\lambda)$ and $K(B_{\ell^n_q})$ whenever $p=1$ and $\lambda=1$.\\[1mm]

In the classical sense, Bombieri and Bourgain \cite{bombieri-1962,bombieri-2004} have studied the constant $K(\mathbb{D},\lambda)$ whenever $\lambda\geq 1$. In fact, Bombieri \cite{bombieri-1962} has obtained the exact value 
\begin{equation*}
	K(\mathbb{D},\lambda)=\frac{1}{3\lambda-2\sqrt{2(\lambda^2-1)}},
\end{equation*}
where $\lambda \in [1,\sqrt{2}]$, and Bombieri and Bourgain \cite{bombieri-2004} have determined the exact asymptotic behavior of $K(\mathbb{D},\lambda)$ as $\sqrt{\lambda^2-1}/\lambda$ whenever $\lambda \rightarrow \infty$.  In 2010, Blasco \cite{Blasco-OTAA-2010} established that $K(\ell^2_q)=0$ for $1\leq q\leq \infty$. Thus, the study of Bohr radius problem becomes irrelevant for $\text{dim}(X)>1$. Hence, it seems reasonable to modify the Bohr inequality to obtain more interesting result. That is why Defant {\it et al.} \cite{defant-adv-math-2012} have introduced $\lambda$ in the definition of Bohr radius $K(B_{\ell^n_{\infty}},T,\lambda)$ of bounded linear operator $T$.\\[1mm]

A complete Reinhardt domain $\Omega \subseteq \mathbb{C}^n$ is a domain in $\mathbb{C}^n$ such that if $z=(z_{1},\ldots,z_{n}) \in \Omega$, then $(\lambda e^{i\theta_{1}}z_{1},\ldots,\lambda e^{i\theta_{n}}z_{n}) \in \Omega$ for all $\lambda \in \overline{\mathbb{D}}$ and all $\theta_{i} \in \mathbb{R}$, $i=1, \ldots,n$. In this  paper, we mostly focus on the complete Reinhardt domain $\Omega=B_{\ell^n _q}$, unit balls in $\ell^n_q$ spaces, where $q \in [1,\infty]$. For $q \in [1,\infty)$,
$$B_{\ell^n _q}= \left\{z \in \mathbb{C}^n:\norm{z}_{q}=\left(\sum_{i=1}^{n}|z_{i}|^q\right)^{1/q}<1\right\},$$ 
and 
$$B_{\ell^n _\infty}= \left\{z \in \mathbb{C}^n: \norm{z}_{\infty}=\sup _{1\leq i \leq n} |z_{i}|<1\right\}.$$
Apart from Bohr radius, there is another interesting variant of multidimensional Bohr radius in $\mathbb{C}^n$, the arithmetic Bohr radius. In 2008, Defant {\it et al.} \cite{defant-QJM-2008,defant-Domain-J angew-2009} introduced this radius and later extensively studied the arithmetic Bohr radius for complete Reinhardt domain in $\mathbb{C}^n$. Defant and his coauthors \cite{defant-Domain-J angew-2009} have used the arithmetic Bohr radius as a very useful technical tool to describe the domains of existence of the monomial expansion of holomorphic functions on a complete Reinhardt domain in a complex Banach space (see also the monograph \cite{defant-book}). It is the right time to introduce the following notion of arithmetic Bohr radius of bounded linear operators between complex Banach spaces which will be the center of discussion of our paper.
\begin{defn}
	Let $1\leq q\leq\infty$, $n\in \mathbb{N}$, and $\mathcal{O}(B_{\ell^n_q},X)$ be the set of all holomorphic mappings on the unit ball $B_{\ell^n_q}$ of $\ell^n_q$ spaces into a complex Banach space $X$. Further, let $1\leq p<\infty$, $\lambda\geq 1$, and $T:X\rightarrow Y$ be a bounded linear operator between complex Banach spaces $X$ and $Y$ such that $\norm{T}<\lambda^{1/p}$. Then the $\lambda_p$-arithmetic Bohr radius of $T$ is defined as
	\begin{equation*}
		A_p(\mathcal{O}(B_{\ell^n_q},X),T,\lambda):= \left\{\frac{1}{n}\sum_{i=1}^{n}r_i \, | \, r\in \mathbb{R}^n_{\geq 0}, \, \forall f \in \mathcal{O}(B_{\ell^n_q},X) : 
		\sum_{\alpha\in \mathbb{N}^n_{0}}\norm{T(x_{\alpha})}^p r^{p\alpha} \leq \lambda \norm{f}^p_{B_{\ell^n_q}} \right\}.
	\end{equation*}
\end{defn}
\noindent
In this paper, we are mainly interested in the following classes of holomorphic functions:\\
\begin{itemize}
	\item $\mathcal{O}(B_{\ell^n_q},X)=\mathcal{H}_{\infty}(B_{\ell^n_q},X)$, the set of all bounded holomorphic functions,\\
	\item $\mathcal{O}(B_{\ell^n_q},X)=\mathcal{P}(^mB_{\ell^n_q},X)$, the set of all $m$-homogeneous polynomials on $\mathbb{C}^n$ as defined on $B_{\ell^n_q}$.
\end{itemize} 

\noindent
We write $A_p(B_{\ell^n_q},T,\lambda)$ for $A_p(\mathcal{H}_{\infty}(B_{\ell^n_q}, X),T,\lambda)$. If $T$ is the identity map on $X$ then we use the notation $A_p(B_{\ell^n_q},X,\lambda)$ and $A_p(B_{\ell^n_q},X)$. For the class $\mathcal{P}(^mB_{\ell^n_q},X)$ we write $A^m_p(B_{\ell^n_q},T,\lambda)$ for $A_p(\mathcal{P}(^mB_{\ell^n_q},X),T,\lambda)$. Defant {\it et al.} \cite{defant-QJM-2008} have studied the arithmetic Bohr radius $A(\Omega)$ for $p=1$ and complete  Reinhardt domain $\Omega$ in $\mathbb{C}^n$. Recently, Kumar and Manna \cite{kumar-2023-arxiv} have studied the arithmetic Bohr radius $A_p(B_{\ell^n_q},T,\lambda)$ for $p=1$. Further, Kumar \cite{kumar-2023} and Allu {\it et al.} \cite{Subhadip-Vasu-P3} have investigated the arithmetic Bohr radius for scalar valued and Banach space valued holomorphic functions respectively. Moreover, arithmetic Bohr radius for holomorphic functions on complete Reinhardt domain in $\mathbb{C}^n$ with positive real part has been investigated by Allu {\it et al.} \cite{Subhadip-Vasu-P4-forum}.\\[1mm]

Our main aim of this paper is to study the arithmetic Bohr radius $A_p(B_{\ell^n_q},T,\lambda)$ of bounded linear operators between complex Banach spaces for any $1\leq p<\infty$ and $\lambda\geq 1$. We establish a close relationship between the arithmetic Bohr radius and the classical Bohr radius. In case of identity operator on infinite dimensional complex Banach space $X$, we give the correct asymptotic behavior of the arithmetic Bohr radius $A_p(B_{\ell^n_q},X,\lambda)$. Further, we give the asymptotic estimates of arithmetic Bohr radius $A_p(B_{\ell^n_q},T,\lambda)$, where $T:\ell_s\hookrightarrow \ell_t$ is an embedding between sequence spaces. As a consequence, we prove the exact asymptotic behavior of arithmetic Bohr radius $A_p(B_{\ell^n_q},T,\lambda)$ for any bounded operator $T:\ell_1 \rightarrow \ell_s$.
\section{Preliminaries}
We employ standard terminologies and notions from classical theory of Banach spaces (see e.g. \cite{lindenstrauss-book-I,lindenstrauss-book-II}). We denote $\mathbb{D}=\{z\in \mathbb{C}:|z|<1\}$ to be the open unit disk in complex plane $\mathbb{C}$. Throughout this paper, we assume that every Banach space $X$ is complex. Let $X^*$ denote the topological dual of Banach space $X$ and we write $B_{X}$ for the open unit ball in $X$. As usual, we denote $\ell_q$ (or $\ell^n_q$), $1\leq q\leq \infty$ for the Banach space of all scalar sequences (or $n$-tuples $z=(z_1, \ldots, z_n)\in \mathbb{C}^n$) endowed with the norm 
\begin{equation*}
	\norm{z}_q:= \left(\sum_{j}|z_j|^q\right)^{\frac{1}{q}}\,\,\, \mbox{for} \,\, 1\leq q<\infty
\end{equation*}
and $\norm{z}_{\infty}:= \sup_{j}|z_j|$ for $q=\infty$. For $x\in X$, the absolute value of $x$ is defined by $|x|=x\vee (-x)$. Recall that a Banach lattice is a Banach space $X$ which is a vector lattice with $|x|\leq |y|$ implies $\norm{x}\leq \norm{y}$ for all $x,y\in X$. A Banach lattice is said to be $r$-concave, where $1\leq r<\infty$, if there is a constant $M>0$ such that 
\begin{equation}\label{Pal-Vasu-P5-pre-e-01}
	\left(\sum_{j=1}^{n}\norm{x_j}^r\right)^{\frac{1}{r}}\leq M \norm{\left(\sum_{j=1}^{n}|x_j|^r\right)^{\frac{1}{r}}}
\end{equation} 
for every arbitrarily chosen finite elements $x_1,\ldots, x_n\in X$. The best such constant satisfying \eqref{Pal-Vasu-P5-pre-e-01} is as usual denoted by $M_r(X)$. For $2\leq r\leq \infty$, a Banach space $X$ is said to have cotype $r$, if there is a constant $C>0$ such that for every choice of finitely many vectors $x_1, \ldots, x_n\in X$, we have
\begin{equation}\label{Pal-Vasu-P5-pre-e-02}
	\left(\sum_{j=1}^{n}\norm{x_j}^r\right)^{\frac{1}{r}} \leq C \left(\int_{0}^{1}\norm{\sum_{j=1}^{n}r_j(t)x_j}^2\, dt\right)^{\frac{1}{2}},
\end{equation}
where $r_j$ is the $j$th Rademacher function on $[0,1]$. Further, the best such constant $C$ satisfying \eqref{Pal-Vasu-P5-pre-e-02} is usually denoted by $C_r(X)$. We write 
\begin{equation*}
	\text{Cot}(X):= \inf \left\{2\leq r\leq \infty: X \,\, \mbox{has cotype}\, r\right\}.
\end{equation*}
For a Banach space $X$ with cotype $\infty$ {\it i.e.}, $\text{Cot}(X)=\infty$, we denote $\frac{1}{\text{Cot}(X)}=0$. We note that the notions of concavity and cotype are closely related concepts in the context of Banach lattices. In fact, a $r$-concave Banach lattice with $r\geq 2$ is of cotype $r$. On the other hand, each Banach lattice of cotype $2$ is $2$-concave, and a Banach lattice of cotype $r>2$ is $s$-concave for all $s>r$.\\

For $m,n\in \mathbb{N}$, we consider the following multi-index set
\begin{equation*}
	\Lambda(m,n)=\left\{\alpha=(\alpha_1,\ldots,\alpha_n)\in \mathbb{N}^n_{0}: |\alpha|=m\right\},
\end{equation*}
where $|\alpha|=\alpha_1+\cdots+\alpha_n$. For $1\leq q\leq \infty$, let $\mathcal{P}(^m\ell^n_q,X)$ be the linear space of all $m$-homogeneous polynomial $Q:\ell^n_q\rightarrow X$ of $n$ complex variables with values in $X$, equipped with uniform (or sup) norm 
\begin{equation*}
	\norm{Q}_{\ell^n_q}:=\sup_{z\in B_{\ell^n_{q}}}\norm{Q(z)}_{X}.
\end{equation*}
Using the above multi-index notations, a polynomial $Q\in \mathcal{P}(^m\ell^n_q,X)$
can be expressed in the form 
\begin{equation*}
	Q(z)=\sum_{\alpha\in\Lambda(m,n)}x_{\alpha}z^{\alpha},
\end{equation*}
where $x_{\alpha}\in X$. If $f:B_{\ell^n_q}\rightarrow X$ is an $X$-valued holomorphic function then we write $f(z)=\sum_{\alpha\in \mathbb{N}^n_{0}}x_{\alpha}(f)z^{\alpha}$ as the monomial series expansion of $f$, where $x_{\alpha}(f)=(\partial^{\alpha}f(0))/\alpha!$ is the $\alpha^{\rm th}$ coefficient of the expansion.\\[1mm]

A bounded linear operator $T:X\rightarrow Y$ between two Banach spaces $X$ and $Y$ is called $(r,s)$-summing, $1\leq r,s<\infty$, if there is a constant $\widetilde{C}>0$ such that for each choice of finitely many vectors $x_1,\ldots,x_n\in X$, we have
\begin{equation}\label{Pal-Vasu-P5-pre-e-03}
	\left(\sum_{j=1}^{n}\norm{T(x_j)}^r\right)^{\frac{1}{r}}\leq \widetilde{C}\sup_{x^*\in B_{X^*}}\left(\sum_{j=1}^{n}|x^*(x_j)|^s\right)^{\frac{1}{s}}.
\end{equation}
Also, the best such constant $\widetilde{C}$ satisfying \eqref{Pal-Vasu-P5-pre-e-03} is denoted by $\pi_{r,s}(T)$. In particular, if $r=s$, we call $T$ is $r$-summing and we write the best constant as $\pi_r(T)$.\\

For two sequences of positive numbers $\{x_n\}$ and $\{y_n\}$, we write $x_n\prec y_n$ if there exists a constant $C>0$ such that $x_n\leq Cy_n$ for every $n\in \mathbb{N}$, and $x_n \asymp y_n$ whenever $x_n \prec y_n$ and $y_n \prec x_n$.
\section{Main Results}
In this paper, we first prove the arithmetic Bohr radius $A_p(B_{\ell^n_q},T,\lambda)$ in terms of\\ $A_p(\cup_{m=1}^{\infty}\mathcal{P}(^mB_{\ell^n_q},X),T,\lambda)$ which shows the connection between arithmetic Bohr radii of two different classes $H_{\infty}(B_{\ell^n_q},X)$ and $\cup_{m=1}^{\infty}\mathcal{P}(^mB_{\ell^n_q},X)$ respectively.
\begin{prop}\label{Pal-Vasu-P6-prop-001}
Let  $1\leq p<\infty$ and $T:X\rightarrow Y$ be a bounded operator between complex Banach spaces $X$ and $Y$ with $\norm{T}<\lambda^{1/p}$. Then we have 
		\begin{enumerate}
		\item[(a)] $\left(\frac{\lambda-\norm{T}^p}{2\lambda-\norm{T}^p}\right)^{\frac{1}{p}}A_p\left(\bigcup^{\infty}_{m=1}\mathcal{P}(^mB_{\ell^n_q},X),T,\lambda\right)\leq A_p(B_{\ell^n_q},T,\lambda)\leq A_p\left(\bigcup^{\infty}_{m=1}\mathcal{P}(^mB_{\ell^n_q},X),T,\lambda\right)$
		\item[(b)] $\left(\frac{\lambda-\norm{T}^p}{\lambda-\norm{T}^p+1}\right)^{\frac{1}{p}}A_p\left(\bigcup^{\infty}_{m=1}\mathcal{P}(^mB_{\ell^n_q},X),T\right)\leq A_p(B_{\ell^n_q},T,\lambda).$
	\end{enumerate}
\end{prop}

As discussed before, our aim is to establish the asymptotic estimates for arithmetic Bohr radius $A_p(B_{\ell^n_q},T,\lambda)$ of bounded linear operators. To prove those estimates, we shall use the bounds of Bohr radius $K^p(B_{\ell^n_q},T,\lambda)$ of bounded linear operators $T$. That is why, we require a relation between these two radii. The next result Lemma \ref{Pal-Vasu-P6-lem-001} will provide a connection between the arithmetic Bohr radius and Bohr radius of bounded linear operators. Before that, we need the following notation from \cite{defant-QJM-2008}. For bounded complete Reinhardt domains $\Omega_1$ and $\Omega_2$ in $\mathbb{C}^n$ we use the notation 
\begin{equation*}
	S(\Omega_1, \Omega_2):= \inf \{t>0: \Omega_1 \subset t\Omega_2 \}.
\end{equation*}
Further, Defant {\it et al.} \cite{defant-2004} have shown that if $\Omega$ is a bounded complete Reinhardt domain in $\mathbb{C}^n$, then 
\begin{align*}
	 &S(\Omega, B_{\ell^n_q})=\sup_{z \in \Omega}\left(\sum_{j=1}^{n}|z_j|^q\right)^{\frac{1}{q}}, \,\,\, 1\leq q<\infty, \\[2mm] & 
	 S(\Omega, B_{\ell^n_{\infty}})= \sup_{z \in \Omega} \sup_{1\leq j\leq n}|z_j|.
\end{align*}
The following lemma is a key tool in estimating the lower bounds for arithmetic Bohr radius of bounded linear operators.
\begin{lem}\label{Pal-Vasu-P6-lem-001}
Let $1\leq p<\infty$ and $T:X\rightarrow Y$ be a bounded linear operator between complex Banach spaces $X$ and $Y$. Then for every $1\leq \lambda <\infty$ and $n\in \mathbb{N}$, we have
\begin{equation*}
	A_p(B_{\ell^n_q},T,\lambda)\geq \frac{K^p(B_{\ell^n_q},T,\lambda)}{n^{1/q}}.
\end{equation*}
\end{lem}
Next, we want to study the arithmetic Bohr radius $A_p(B_{\ell^n_q},T,\lambda)$ for identity operator $T$ on infinite dimensional complex Banach space $X$. Recently, Allu {\it et al.} \cite{Subhadip-Vasu-P3,Subhadip-Vasu-P5} have studied the correct asymptotic estimates for $A_p(B_{\ell^n_q},X,\lambda)$ for finite dimensional complex Banach space $X$. To establish the estimates for arithmetic Bohr radius $A_p(B_{\ell^n_q},X,\lambda)$, we need the following result for Bohr radius $K^p(B_{\ell^n_q},X,\lambda)$, where $X$ is an infinite dimensional complex Banach space.
\begin{customthm}{A}\label{Pal-Vasu-P5-thm-A}\cite{Subhadip-Vasu-P5}
	Let $1\leq p<\infty$ and $1\leq q\leq \infty$. If $X$ is an infinite dimensional complex Banach space of cotype $r$ with $2\leq r\leq \infty$, then we have 
	$$
	K^p(B_{\ell^n_q},X,\lambda)\leq \left\{\begin{array}{ll}
		\frac{\lambda^{1/p}}{n^{\frac{1}{p}-\frac{1}{q}}}, & \mbox{for \, $q\leq \text{Cot}(X)=r$}\\[5mm] 
		
		\frac{\lambda^{1/p}}{n^{\frac{1}{p}\left(1-\frac{1}{Cot(X)}\right)-\frac{1}{q}\left(1-\frac{1}{p}\right)}}, & \mbox{for \, $q>\text{Cot}(X)=r$}
	\end{array}\right.
	$$
	and 
	$$
	K^p(B_{\ell^n_q},X,\lambda)\geq \left\{\begin{array}{ll}
		\left(\frac{\lambda-1}{\lambda}\right)^{1/p}\frac{1}{e\, C_r(X)}, & \mbox{for \, $q\leq s$}\\[5mm] 
		
		\left(\frac{\lambda-1}{\lambda}\right)^{1/p}\frac{1}{e\,C_r(X)n^{\frac{1}{s}-\frac{1}{q}}}, & \mbox{for \, $q>s$},
	\end{array}\right.
	$$
	where $s=r/(r-1)$.
\end{customthm}
We observe that asymptotic estimates of arithmetic Bohr radius $A_p(B_{\ell^n_q},X,\lambda)$ is governed by the concept of geometry of Banach space. In fact, the bounds of arithmetic Bohr radius are influenced by optimal cotype $\text{Cot}(X):=\mbox{the infimum over all $r$, where}\, \, r\in [2,\infty)$ for which the Banach space $X$ has cotype $r$. 
\begin{thm}\label{Pal-Vasu-P6-thm-001}
	Let $1\leq p<\infty$ and $1\leq \lambda <\infty$. Suppose $X$ is an infinite dimensional complex Banach space of cotype $r$ with $2\leq r\leq \infty$. Then for any $1\leq q\leq \infty$, we have 
		$$
	A_p(B_{\ell^n_q},X,\lambda)\leq \left\{\begin{array}{ll}
		\frac{\lambda^{1/p}}{n}, & \mbox{for \, $p\leq r$}\\[5mm] 
		
	\frac{\lambda^{1/p}}{n^{1-p\left(\frac{1}{\text{Cot}(X)}-\frac{1}{q}\right)}}, & \mbox{for \, $p>r$},
	\end{array}\right.
	$$
	and 
		$$
	A_p(B_{\ell^n_q},X,\lambda)\geq \left\{\begin{array}{ll}
	\left(\frac{\lambda-1}{\lambda}\right)^{1/p}\frac{1}{e\, C_r(X)n^{1/q}}, & \mbox{for \, $q\leq s$}\\[5mm] 
		
	\left(\frac{\lambda-1}{\lambda}\right)^{1/p}\frac{1}{e\, C_r(X)n^{1/s}}, & \mbox{for \, $q>s$},
	\end{array}\right.
	$$
	where $s=r/(r-1)$.
\end{thm}
Following the techniques from Defant {\it et al.} \cite{defant-adv-math-2012} we investigate the Bohr radii of operators between Banach spaces, specifically here we consider the embeddings $T:\ell_s \hookrightarrow \ell_t, \,\, 1\leq s\leq t<\infty$ and arbitrary operators $T: \ell_1 \rightarrow \ell_s$. Our first aim is to study the Bohr radius $K^p(B_{\ell^n_q},T,\lambda)$ for aforesaid operators. Before proceeding  further discussion, we want to mention an important result recently proved in \cite{Subhadip-Vasu-P5} estimating the lower bounds of $K^p(B_{\ell^n_q},T,\lambda)$ for any bounded linear operator $T$ and $1\leq p<\infty$ which plays a integral role in the following results.
\begin{customthm}{B}\label{Pal-Vasu-P5-thm-B}\cite{Subhadip-Vasu-P5}
		Let $1\leq p<\infty$ and $T:X\rightarrow Y$ be a bounded operator between Banach spaces $X$ and $Y$.
	\begin{enumerate}
		\item Assume that $X$ or $Y$ is of cotype $r$ with $2\leq r\leq \infty$. Then there is a constant $D>0$ such that for every $\norm{T}<\lambda^{1/p}$, $1\leq q\leq \infty$, and $n$
		$$
		K^p(B_{\ell^n_q},T, \lambda) \geq \left\{\begin{array}{ll}
			D\, \left(\frac{\lambda - \norm{T}^p}{\lambda}\right)^{\frac{1}{p}}, & \mbox{for \, $q\leq s$}\\[5mm] 
			
			D\, \left(\frac{\lambda - \norm{T}^p}{\lambda}\right)^{\frac{1}{p}} \left(\frac{1}{n}\right)^{\frac{1}{s}-\frac{1}{q}}, & \mbox{for \, $q>s$},
		\end{array}\right.
		$$
		where $s=r/(r-1)$.
		\item Assume that $Y$ is a $r$-concave Banach lattice with $2\leq r<\infty$ and there is a $1\leq s<r$ such that the operator $T$ is $(s,1)$-summing. Then there is a constant $D>0$  such that 
		\begin{equation*}
			K^p(B_{\ell^n_q},T, \lambda) \geq D\, \left(\frac{\lambda-\norm{T}^p}{2\lambda -\norm{T}^p}\right)^{\frac{1}{p}} \left(\frac{\log n}{n}\right)^{1-\frac{1}{r}}
		\end{equation*}
		for every $\norm{T}<\lambda^{1/p}$, $1\leq q\leq \infty$, and $n\in \mathbb{N}$.
	\end{enumerate}
\end{customthm}
By unifying Theorem \ref{Pal-Vasu-P5-thm-B} and the well-known Grothendieck type theorems, we execute the proof of Theorem \ref{Pal-Vasu-P6-thm-002} and Theorem \ref{Pal-Vasu-P6-thm-003}. The actual form of Grothendieck's theorem \cite{diestel-Abs-sum-book} states that every bounded linear operator $T:\ell_1 \rightarrow \ell_2$ is absolutely summing {\it i.e.,} $1$-summing. A preceding version of this result is actually an old result by Littlewood \cite{littlewood-1930}: The formal inclusion $\ell_1 \hookrightarrow \ell_{4/3}$ is $(4/3, 1)$-summing, also known as the Littlewood's $4/3$-theorem. After this, these classical results got many extensions in $\ell_q$ spaces, and the most significant ones are due to Kwapie\'{n} \cite{kwapien-1968} and Bennett-Carl \cite{Bennett-1973,Carl-1974}. In Theorem \ref{Pal-Vasu-P6-thm-002} and Theorem \ref{Pal-Vasu-P6-thm-003}, we obtain the estimates for arithmetic Bohr radius of inclusions $T: \ell_s \hookrightarrow\ell_t$ and bounded operators $T:\ell_1 \rightarrow \ell_s$ respectively.

\begin{thm}\label{Pal-Vasu-P6-thm-002}
	Let $1\leq p<\infty$, $1\leq \lambda<\infty$, and $1\leq s<t<\infty$. Then with constants depending only on $\lambda$, $p$, and $s,t$, we have
	$$
	K^p(B_{\ell^n_q},\ell_s\hookrightarrow\ell_t,\lambda)\prec \left\{\begin{array}{ll}
		n^{1-\frac{1}{p}}\left(\frac{\log n}{n}\right)^{1-\frac{1}{\min\{q,2\}}}, & \mbox{for \, $1\leq s\leq 2$}\\[5mm] 
		
		\frac{1}{n^{\frac{1}{p}-\frac{1}{q}}}, & \mbox{for \, $q\leq s$\, and\, $2\leq s<\infty$}\\[5mm]
		 	\frac{1}{n^{\frac{1}{p}-\frac{1}{s}}}, & \mbox{for \, $q> s$\, and\, $2\leq s<\infty$}
	\end{array}\right.
	$$
and 
 $$
	 K^p(B_{\ell^n_q},\ell_s\hookrightarrow\ell_t,\lambda)\succ \left\{\begin{array}{ll}
	 	\sqrt{\frac{\log n}{n}}, & \mbox{for \, $1\leq s\leq 2$}\\[5mm] 
	 	
	 	\frac{1}{e}, & \mbox{for \, $q\leq r$\, and\, $2\leq s<\infty$}\\[5mm]
	 	\frac{1}{en^{\frac{1}{r}-\frac{1}{q}}}, & \mbox{for \, $q> r$\, and\, $2\leq s<\infty$},
	 \end{array}\right.
	 $$
	 where $r=s/(s-1)$.
\end{thm}

\begin{thm}\label{Pal-Vasu-P6-thm-003}
	Let $1\leq p<\infty$, $1\leq s<\infty$, and $T:\ell_{1} \rightarrow \ell_{s}$ be any bounded operator. Then with constants depending only on $T$, $\lambda$ and $p$, we have
	\begin{equation*}
		\left(\frac{\log n}{n}\right)^{1-\frac{1}{\max\{2,s\}}} \prec K^p(B_{\ell^n_q},T,\lambda) \prec n^{1-\frac{1}{p}}\left(\frac{\log n}{n}\right)^{1-\frac{1}{\min\{q,2\}}}.
	\end{equation*}
\end{thm}
We note that the classical multidimensional Bohr radius and arithmetic Bohr radius are two closely related variants of Bohr radius. To investigate the asymptotic behavior of arithmetic Bohr radius, we find the estimates of classical Bohr radius very useful. We use this observation for our next main results too. By virtue of Theorem \ref{Pal-Vasu-P6-thm-002} and Theorem \ref{Pal-Vasu-P6-thm-003}, we study the asymptotic estimates for arithmetic Bohr radius $A_p(B_{\ell^n_q},T,\lambda)$, where $T$ is the inclusions $\ell_s \hookrightarrow \ell_t, \,\, 1\leq s<t<\infty$ and any bounded operators $\ell_1 \rightarrow \ell_s, \,\, 1\leq s<\infty$.

\begin{thm}\label{Pal-Vasu-P6-thm-004}
		Let $1\leq p<\infty$, $1\leq \lambda<\infty$, and $1\leq s<t<\infty$. Then with constants depending only on $\lambda$, $p$, and $s,t$, we have
	$$
	A_p(B_{\ell^n_q},\ell_s\hookrightarrow\ell_t,\lambda)\prec \left\{\begin{array}{ll}
		\frac{(\log n)^{1-\frac{1}{\min\{q,2\}}}}{n^{\frac{1}{p}+\frac{1}{\max\{q,2\}}-\frac{1}{2}}}, & \mbox{for \, $1\leq s\leq 2$}\\[5mm] 
		
	\frac{1}{n^{1-\frac{p}{q}}}, & \mbox{for \, $q\leq s$\, and\, $2\leq s<\infty$}\\[5mm]
	\frac{1}{n^{1-\frac{p}{s}}}, & \mbox{for \, $q> s$\, and\, $2\leq s<\infty$}
	\end{array}\right.
	$$
	and 
	$$
	A_p(B_{\ell^n_q},\ell_s\hookrightarrow\ell_t,\lambda)\succ \left\{\begin{array}{ll}
	\frac{\sqrt{\log n}}{n^{\frac{1}{2}+\frac{1}{q}}}, & \mbox{for \, $1\leq s\leq 2$}\\[5mm] 
		
		\frac{1}{en^{\frac{1}{q}}}, & \mbox{for \, $q\leq r$\, and\, $2\leq s<\infty$}\\[5mm]
		\frac{1}{en^{\frac{1}{r}-\frac{1}{q}}}, & \mbox{for \, $q> r$\, and\, $2\leq s<\infty$},
	\end{array}\right.
	$$
	where $r=s/(s-1)$.
\end{thm}

\begin{thm}\label{Pal-Vasu-P6-thm-005}
	Let $1\leq p<\infty$, $1\leq s<\infty$, and $T:\ell_{1} \rightarrow \ell_{s}$ be any bounded operator. Then with constants depending only on $T$, $\lambda$ and $p$, we have
	\begin{equation*}
		\frac{1}{n^{1/q}}\left(\frac{\log n}{n}\right)^{1-\frac{1}{\max\{2,s\}}} \prec A_p(B_{\ell^n_q},T,\lambda) \prec 	\frac{(\log n)^{1-\frac{1}{\min\{q,2\}}}}{n^{\frac{1}{p}+\frac{1}{\max\{q,2\}}-\frac{1}{2}}}.
	\end{equation*}
\end{thm}

\section{Proof of the Main results}
\begin{pf}[{\bf Proof of Proposition \ref{Pal-Vasu-P6-prop-001}}]
	The right-hand inequality of (a) is clear from the fact that 
	\begin{equation*}
		\bigcup^{\infty}_{m=1}\mathcal{P}(^mB_{\ell^n_q},X)\subset \mathcal{H}_{\infty}(B_{\ell^n_q},X).
	\end{equation*}
	Therefore, we want to prove the left-hand inequality of (a). Let $r=(r_1,\ldots,r_n)\in \mathbb{R}^n_{\geq 0}$ be such that 
	\begin{equation*}
		\sum_{\alpha\in \Lambda(m,n)}\norm{T(x_{\alpha}(Q))}^pr^{p\alpha}\leq \lambda \sup_{z\in B_{\ell^n_q}}\norm{Q(z)}^p
	\end{equation*}
	for all vector-valued $m$-homogeneous polynomial $Q(z)=\sum_{\alpha\in \Lambda(m,n)}x_{\alpha}(Q)z^{\alpha}\in \mathcal{P}(^mB_{\ell^n_q},X)$.
	Fix some $f(z)=\sum_{\alpha\in \mathbb{N}^n_{0}}x_{\alpha}z^{\alpha} \in \mathcal{H}_{\infty}(B_{\ell^n_q},X)$. It is worth noting that for $1\leq p<\infty$,
	\begin{equation*}
		0<	\left(\frac{\lambda-\norm{T}^p}{2\lambda-\norm{T}^p}\right)^{1/p}<1.
	\end{equation*}
	Then for $1\leq p<\infty$ and by using the Cauchy-Riemann inequalities, we obtain
	\begin{align*}
		&\sum_{\alpha\in \mathbb{N}^n_{0}}\norm{T(x_{\alpha})}^p \left(\frac{\lambda-\norm{T}^p}{2\lambda-\norm{T}^p}\right)r^{p\alpha}\\[2mm]&= \norm{T(x_{0})}^p + \sum_{m=1}^{\infty}\left(\frac{\lambda-\norm{T}^p}{2\lambda-\norm{T}^p}\right)^m \sum_{\alpha\in \Lambda(m,n)}\norm{T(x_{\alpha})}^p r^{p\alpha}\\[2mm] & \leq 
		\norm{T}^p\norm{x_0}^p + \sum_{m=1}^{\infty}\left(\frac{\lambda-\norm{T}^p}{2\lambda-\norm{T}^p}\right)^m \lambda \norm{\sum_{\alpha\in \Lambda(m,n)}x_{\alpha}z^{\alpha}}^p_{B_{\ell^n_q}}\\[2mm] & \leq 
		\left(\norm{T}^p+\lambda \sum_{m=1}^{\infty}\left(\frac{\lambda-\norm{T}^p}{2\lambda-\norm{T}^p}\right)^m\right)\norm{f}^p_{B_{\ell^n_q}}=\lambda \norm{f}^p_{B_{\ell^n_q}}.
	\end{align*}
	Hence, it follows that
	\begin{equation*}
		\left(\frac{\lambda-\norm{T}^p}{2\lambda-\norm{T}^p}\right) \frac{1}{n}\sum_{i=1}^{n}r_i\leq A_p(B_{\ell^n_q},T,\lambda),
	\end{equation*}
	which gives the desired result.\\[2mm]
	
	By following the similar process, we shall now prove the inequality (b). Let $r=(r_1,\ldots,r_n)\in \mathbb{R}^n_{\geq 0}$ be such that for all $1\leq p<\infty$ and $m\in \mathbb{N}$,
	\begin{equation*}
		\sum_{\alpha\in \Lambda(m,n)}\norm{T(x_{\alpha}(Q))}^pr^{p\alpha}\leq \sup_{z\in B_{\ell^n_q}} \norm{Q(z)}^p
	\end{equation*}
	for all $Q(z)=\sum_{\alpha\in \mathbb{N}^n_{0}}x_{\alpha}z^{\alpha}\in \mathcal{P}(^mB_{\ell^n_q},X)$. Suppose $g(z)=\sum_{\alpha\in \mathbb{N}^n_{0}}x_{\alpha}z^{\alpha}$ is a bounded holomorphic function on $B_{\ell^n_q}$ with values in complex Banach space $X$. By the virtue of Cauchy-Riemann inequalities, we have
	\begin{align*}
		& \sum_{\alpha\in \mathbb{N}^n_{0}}\left(\frac{\lambda-\norm{T}^p}{\lambda-\norm{T}^p+1}\right)^{\alpha}r^{p\alpha} \\[2mm]& = \norm{T(x_0)}^p + \sum_{m=1}^{\infty} \left(\frac{\lambda-\norm{T}^p}{\lambda-\norm{T}^p+1}\right)^m \sum_{\alpha\in \Lambda(m,n)} \norm{T(x_{\alpha})}^pr^{p\alpha}\\[2mm]& \leq 
		\norm{T}^p\norm{x_0}^p + \sum_{m=1}^{\infty} \left(\frac{\lambda-\norm{T}^p}{\lambda-\norm{T}^p+1}\right)^m \sup_{z\in B_{\ell^n_q}} \norm{g(z)}^p\\[2mm]& \leq
		\left(\norm{T}^p+\sum_{m=1}^{\infty}\left(\frac{\lambda-\norm{T}^p}{\lambda-\norm{T}^p+1}\right)^m\right)\norm{g}^p_{B_{\ell^n_q}}=\lambda \norm{g}^p_{B_{\ell^n_q}}.
	\end{align*}
	Evidently, we obtain 
	\begin{equation*}
		\left(\frac{\lambda-\norm{T}^p}{\lambda-\norm{T}^p+1}\right)^{1/p}\frac{1}{n}\sum_{i=1}^{n}r_i\leq A_p(B_{\ell^n_q},T,\lambda),
	\end{equation*}
	which gives the required inequality. This completes the proof. 
\end{pf}

\begin{pf}[{\bf Proof of Lemma \ref{Pal-Vasu-P6-lem-001}}]
	From the definition of the notation $S$, we have 
	\begin{equation*}
		S(B_{\ell^n_q}, B_{\ell^n_1})=\sup_{z\in B_{\ell^n_q}} \norm{z}_{\ell^n_1}=n^{1-\frac{1}{q}}.
	\end{equation*}
	Then for a given $0<\epsilon<K^p(B_{\ell^n_q},T,\lambda)$ we can find a member $\widetilde{z_{0}}=(\widetilde{z_1}, \ldots, \widetilde{z_n})\in B_{\ell^n_q}$ such that 
	\begin{equation*}
		\norm{\widetilde{z_{0}}}_{\ell^n_1}\geq S(B_{\ell^n_q}, B_{\ell^n_1})-\epsilon.
	\end{equation*}
	Let $s:= K^p(B_{\ell^n_q},T,\lambda)-\epsilon$, $w=s\widetilde{z_{0}}$, and $r:=s|\widetilde{z_{0}}|=|w|=s(|\widetilde{z_1}|,\ldots, |\widetilde{z_n}|)$. Since $w\in sB_{\ell^n_q}$ and $s<K^p(B_{\ell^n_q},T,\lambda)$, for $f(z)=\sum_{\alpha\in \mathbb{N}^n_{0}}x_{\alpha}z^{\alpha}\in \mathcal{H}_{\infty}(B_{\ell^n_q},X)$ we obtain
	\begin{equation*}
		\sum_{\alpha\in \mathbb{N}^n_{0}}\norm{T(x_{\alpha})}^p r^{p\alpha} = \sum_{\alpha\in \mathbb{N}^n_{0}} \norm{T(x_{\alpha})s^{|\alpha|}\widetilde{z_{0}}}^p \leq \norm{\sum_{\alpha\in \mathbb{N}^n_{0}}\norm{T(x_{\alpha})z^{\alpha}}}^p_{sB_{\ell^n_q}}\leq \lambda \norm{f}^p_{B_{\ell^n_q}}.
	\end{equation*}
	It follows that 
	\begin{equation*}
		\frac{1}{n}\sum_{j=1}^{n}r_j=\frac{K^p(B_{\ell^n_q},T,\lambda)-\epsilon}{n} \norm{\widetilde{z_{0}}}_{\ell^n_1}\geq \frac{K^p(B_{\ell^n_q},T,\lambda)-\epsilon}{n} \left(n^{1-\frac{1}{q}}-\epsilon\right),
	\end{equation*}
	which holds for all $\epsilon>0$.
	Therefore, we obtain 
	\begin{equation*}
		A_p(B_{\ell^n_q},T,\lambda) \geq \frac{K^p(B_{\ell^n_q},T,\lambda)}{n^{1/q}}.
	\end{equation*}
	This completes the proof.
\end{pf}

\begin{pf}[{\bf Proof of Theorem \ref{Pal-Vasu-P6-thm-001}}]
	The lower bound of $A_p(B_{\ell^n_q},X,\lambda)$ is an immediate consequence of Lemma \ref{Pal-Vasu-P6-lem-001} and Theorem \ref{Pal-Vasu-P5-thm-A}. Indeed, the authors \cite{Subhadip-Vasu-P5} have shown the following: If $X$ is an infinite dimensional complex Banach space of cotype $r$ with $2\leq r\leq \infty$ then for every $1\leq p<\infty$ we have 
	$$
	K^p(B_{\ell^n_q},X,\lambda)\geq \left\{\begin{array}{ll}
		\left(\frac{\lambda-1}{\lambda}\right)^{1/p}\frac{1}{e\, C_r(X)}, & \mbox{for \, $q\leq s$}\\[5mm] 
		
		\left(\frac{\lambda-1}{\lambda}\right)^{1/p}\frac{1}{e\,C_r(X)n^{\frac{1}{s}-\frac{1}{q}}}, & \mbox{for \, $q>s$},
	\end{array}\right.
	$$
	where $s=r/(r-1)$. In view of the above lower estimate of $K^p(B_{\ell^n_q},X,\lambda)$ and Lemma \ref{Pal-Vasu-P6-lem-001}, we obtain
	$$
	A_p(B_{\ell^n_q},X,\lambda)\geq \left\{\begin{array}{ll}
		\left(\frac{\lambda-1}{\lambda}\right)^{1/p}\frac{1}{e\, C_r(X)n^{1/q}}, & \mbox{for \, $q\leq s$}\\[5mm] 
		
		\left(\frac{\lambda-1}{\lambda}\right)^{1/p}\frac{1}{e\, C_r(X)n^{1/s}}, & \mbox{for \, $q>s$},
	\end{array}\right.
	$$
	where $s=r/(r-1)$.\\[2mm]
	\noindent
	To prove the upper bound we require a deep result by Maurey and Pisier \cite[Theorem 14.5]{diestel-Abs-sum-book}: For every infinite dimensional Banach space $X$,
	\begin{equation*}
		\inf \left\{2\leq r\leq \infty : X \,\, \mbox{has cotype}\,\, r \right\}=\sup \left\{2\leq r\leq \infty : X \,\,\mbox{finitely factors}\,\, \ell_r \hookrightarrow \ell_{\infty}\right\}.
	\end{equation*}
	We say that a Banach space $X$ finitely factors $\ell_r \hookrightarrow \ell_{\infty}$ if for every $\epsilon>0$ and $n\in \mathbb{N}$, there exist $x_1, \ldots, x_n\in X$ such that 
	\begin{equation}\label{Pal-Vasu-P6-eqn-001}
		\frac{1}{1+\epsilon}\norm{z}_{\infty}\leq \norm{\sum_{j=1}^{n}z_jx_j}\leq \norm{z}_r
	\end{equation}
	for all $z=(z_1, \ldots, z_n)\in \mathbb{C}^n$. In particular, by considering $z=e_j=(0,\ldots, 0,1,0,\ldots,0)$, \eqref{Pal-Vasu-P6-eqn-001} becomes
	\begin{equation}\label{Pal-Vasu-P6-eqn-002}
		\frac{1}{1+\epsilon}=\frac{1}{1+\epsilon}\norm{e_j}_{\infty}\leq \norm{x_j},\,\,\, 1\leq j\leq n.
	\end{equation}
	Let $r=(r_1, \ldots,r_n)\in \mathbb{R}^n_{\geq 0}$ be such that 
	\begin{equation*}
		\sum_{j=1}^{\infty}\norm{x_j(Q)}^p r^p_j \leq \sup_{z\in B_{\ell^n_q}} \norm{Q(z)}^p
	\end{equation*}
	for any $1\leq p<\infty$ and for every $1$-homogeneous polynomial $Q(z)=\sum_{j}x_j(Q)z_j \in \mathcal{P}(^1B_{\ell^n_q},X)$, where $X$ is an infinite dimensional complex Banach space of cotype $r$ with $2\leq r\leq \infty$. Therefore, in view of \eqref{Pal-Vasu-P6-eqn-002}, for a given $\epsilon>0$, there exist $x_1, \ldots, x_n \in X$ such that for any $1\leq p<\infty$ we have 
	\begin{equation}\label{Pal-Vasu_P6-eqn-003}
		\frac{1}{(1+\epsilon)^p}\sum_{j=1}^{n} r_j \leq \sum_{j=1}^{n}\norm{x_j}^p r^p_j,
	\end{equation}
	where $r_j \in \mathbb{R}^n_{\geq 0}$, $1\leq j\leq n$. Further, the above equation \eqref{Pal-Vasu_P6-eqn-003} gives 
	\begin{equation}\label{Pal-Vasu_P6-eqn-004}
		\frac{n}{(1+\epsilon)^p} \left(\frac{1}{n}\sum_{j=1}^{\infty}r_j\right) \leq \sum_{j=1}^{n}\norm{x_j}^p r^p_j \leq \sup_{z\in B_{\ell^n_q}} \norm{\sum_{j=1}^{n}z_j x_j}^p \leq \sup_{z\in B_{\ell^n_q}} \norm{z}^p_r.
	\end{equation}
	Now, we observe that if $q\leq \text{Cot}(X)$ then $\sup_{z\in B_{\ell^n_q}}\norm{z}_{\text{Cot}(X)}\leq 1$. Then, by the virtue of \eqref{Pal-Vasu_P6-eqn-004}, we obtain
	\begin{equation*}
		\frac{1}{n}\sum_{j=1}^{n}r_j \leq \frac{(1+\epsilon)^p}{n},
	\end{equation*}
	which holds for every $\epsilon>0 $ and $1\leq p<\infty$. Hence, it shows that 
	\begin{equation}\label{Pal-Vasu-P6-eqn-014}
		A^1_p(B_{\ell^n_q},X)\leq \frac{1}{n}.
	\end{equation}
	On the other hand, if $q>\text{Cot}(X)=r$, then we obtain
	\begin{equation*}
		\norm{z}_{\text{Cot}(X)}\leq n^{\frac{1}{\text{Cot}(X)}-\frac{1}{q}}\norm{z}_q. \quad (\mbox{ we assume the convention}\,\, n^{\frac{1}{\infty}}=1)
	\end{equation*}
	Therefore, in view of \eqref{Pal-Vasu_P6-eqn-004} we have the following:
	\begin{equation*}
		\frac{1}{n}\sum_{j=1}^{n}r_j \leq \frac{(1+\epsilon)^p}{n^{1-p\left(\frac{1}{\text{Cot}(X)}-\frac{1}{q}\right)}}
	\end{equation*}
	holds for every $\epsilon>0$. Consequently, it gives that
	\begin{equation}\label{Pal-Vasu-P6-eqn-015}
		A^1_p(B_{\ell^n_q},X)\leq \frac{1}{n^{1-p\left(\frac{1}{\text{Cot}(X)}-\frac{1}{q}\right)}}.
	\end{equation}
Now, it is worth noting that 
	\begin{equation}\label{Pal-Vasu-P6-eqn-016}
		A_p(B_{\ell^n_q},X,\lambda)\leq A^1_p(B_{\ell^n_q},X,\lambda)=\lambda^{1/p}A_p(B_{\ell^n_q},X).
	\end{equation}
Therefore, by using \eqref{Pal-Vasu-P6-eqn-014} and \eqref{Pal-Vasu-P6-eqn-015} in \eqref{Pal-Vasu-P6-eqn-016}, we obtain
	$$
A_p(B_{\ell^n_q},X,\lambda)\leq \left\{\begin{array}{ll}
	\frac{\lambda^{1/p}}{n}, & \mbox{for \, $p\leq r$}\\[5mm] 
	
	\frac{\lambda^{1/p}}{n^{1-p\left(\frac{1}{\text{Cot}(X)}-\frac{1}{q}\right)}}, & \mbox{for \, $p>r$}.
\end{array}\right.
$$
This completes the proof.
\end{pf}

\begin{pf}[{\bf Proof of Theorem \ref{Pal-Vasu-P6-thm-002}}]
	We first establish the upper bound of $K^p(B_{\ell^n_q},\ell_s\hookrightarrow\ell_t,\lambda)$. We divide the proof into several cases.\\[1mm]
	{\it The case $1\leq s\leq 2$.} To establish the upper estimate in this case, we shall use the following observation:
	\begin{equation}\label{Pal-Vasu-P6-eqn-008}
		K^p(B_{\ell^n_q},T,\lambda)\leq K^p\left(B_{\ell^n_q},\mathbb{C},\frac{\lambda}{\norm{T}^p}\right),
	\end{equation}
	which holds for any bounded operator $T$ and for each $p\in [1,\infty)$. Moreover, in \cite[Theorem 3.1]{Subhadip-Vasu-P5}, the authors have proved that 
	\begin{equation}\label{Pal-Vasu-P6-eqn-008-a}
		K^p(B_{\ell^n_q},X,\lambda)\prec n^{1-\frac{1}{p}}\left(\frac{\log n}{n}\right)^{1-\frac{1}{\min\{q,2\}}}
	\end{equation}
	for every finite dimensional complex Banach space $X$ and for each $1\leq p<\infty$.
	Hence, in view of \eqref{Pal-Vasu-P6-eqn-008} and \eqref{Pal-Vasu-P6-eqn-008-a}, we obtain the desired upper estimate.\\[1mm]
	{\it The case $2\leq s<\infty$.} By the virtue of result by Maurey and Pisier \cite[Theorem 14.5]{diestel-Abs-sum-book} and by use of the definition of $K^p_{1}(B_{\ell^n_q},\ell_s\hookrightarrow\ell_t)$ we have
	\begin{align}\label{Pal-Vasu-P6-eqn-005}
		n^{1-\frac{p}{q}}=\sum_{j=1}^{n}\norm{e_j}^p_{s}\frac{1}{n^{\frac{p}{q}}}\nonumber&
		\leq \frac{1}{(K^p_{1}(B_{\ell^n_q},\ell_s\hookrightarrow \ell_t))^p} \sup_{z\in B_{\ell^n_q}}\norm{\sum_{j=1}^{n}z_j e_j}^p_{s}\\[2mm]&
		\leq \frac{1}{(K^p_{1}(B_{\ell^n_q},\ell_s\hookrightarrow \ell_t))^p} \sup_{z\in B_{\ell^n_q}} \norm{z}^p_s.
	\end{align}
	Now, we observe that if $q\leq s$ then $\sup_{z\in B_{\ell^n_q}}\norm{z}_s\leq 1$. Thus, in view of \eqref{Pal-Vasu-P6-eqn-005}, we have 
	\begin{equation}\label{Pal-Vasu-P6-eqn-006}
		K^p_{1}(B_{\ell^n_q},\ell_s\hookrightarrow\ell_t) \leq \frac{1}{n^{\frac{1}{p}-\frac{1}{q}}}.
	\end{equation}
	Conversely, if $q>s$ we have 
	\begin{equation*}
		\norm{z}_s\leq \norm{z}_q n^{\frac{1}{s}-\frac{1}{q}},
	\end{equation*}
	whence, from \eqref{Pal-Vasu-P6-eqn-005} we obtain
	\begin{equation}\label{Pal-Vasu-P6-eqn-007}
		K^p_{1}(B_{\ell^n_q},\ell_s\hookrightarrow\ell_t) \leq \frac{1}{n^{\frac{1}{p}-\frac{1}{s}}}.
	\end{equation}
	We note that $K^p_{1}(B_{\ell^n_q},T,\lambda)=\lambda^{1/p}K^p_{1}(B_{\ell^n_q},T)$ for any bounded operator $T$. Therefore, in view of \eqref{Pal-Vasu-P6-eqn-006} and \eqref{Pal-Vasu-P6-eqn-007}, we obtain
	$$
	K^p_{1}(B_{\ell^n_q},\ell_s\hookrightarrow\ell_t,\lambda)\leq \left\{\begin{array}{ll}
		\frac{\lambda^{1/p}}{n^{\frac{1}{p}-\frac{1}{q}}}, & \mbox{for \, $q\leq s$}\\[5mm] 
		
		\frac{\lambda^{1/p}}{n^{\frac{1}{p}-\frac{1}{s}}}, & \mbox{for \, $q>s$}.
	\end{array}\right.
	$$	
	We note that $K^p(B_{\ell^n_q},\ell_s\hookrightarrow \ell_t, \lambda)\leq K^p_1(B_{\ell^n_q}, \ell_s \hookrightarrow \ell_t, \lambda)$.
	Therefore, we obtain
	$$
	K^p(B_{\ell^n_q},\ell_s\hookrightarrow\ell_t,\lambda)\prec \left\{\begin{array}{ll}
		\frac{1}{n^{\frac{1}{p}-\frac{1}{q}}}, & \mbox{for \, $q\leq s$\, and\, $2\leq s<\infty$}\\[5mm] 
		
		\frac{1}{n^{\frac{1}{p}-\frac{1}{s}}}, & \mbox{for \, $q>s$\, and\, $2\leq s<\infty$}.
	\end{array}\right.
	$$
	\vspace{2mm}
	
	Next, we proceed to prove the lower bound of $K^p(B_{\ell^n_q},\ell_s\hookrightarrow\ell_t,\lambda)$. To determine the lower bound we consider three different cases.\\[1mm]
	{\it The case $1\leq s<t\leq 2$.} In this case, we use the Bennett-Carl theorem \cite{Bennett-1973,Carl-1974} which states that the inclusion $\ell_s \hookrightarrow \ell_t$ is $(h,1)$-summing where 
	\begin{equation*}
		\frac{1}{h}=\frac{1}{2}+\frac{1}{s}-\max\left\{\frac{1}{t},\frac{1}{2}\right\}.
	\end{equation*}
	It is worth noting that $\ell_t$ is $2$-concave and thus the lower estiamte for $K^p(B_{\ell^n_q},\ell_s\hookrightarrow\ell_t,\lambda)$ is a consequence of Theorem \ref{Pal-Vasu-P5-thm-B}(2). Therefore, we obtain
	\begin{equation*}
		K^p(B_{\ell^n_q},\ell_s\hookrightarrow\ell_t,\lambda)\succ \sqrt{\frac{\log n}{n}}.
	\end{equation*}
	{\it The case $1\leq s<2\leq t$.} Here, we observe that $K^p(B_{\ell^n_q},\ell_s\hookrightarrow\ell_t,\lambda)\geq K^p(B_{\ell^n_q},\ell_s\hookrightarrow\ell_2,\lambda)$ for any $1\leq p<\infty$. Therefore, we obtain the required estimate using the preceding case.\\[1mm]
	{\it The case $2\leq s$.} We notice the following observation:
	\begin{equation*}
		K^p(B_{\ell^n_q},\ell_s\hookrightarrow\ell_t,\lambda)\geq K^p(B_{\ell^n_q},\ell_s\hookrightarrow\ell_s,\lambda).
	\end{equation*}
	Then, in view of Theorem \ref{Pal-Vasu-P5-thm-A}, we obtain the concluding part of our required estimate. This completes the proof.
\end{pf}

\begin{pf}[{\bf Proof of Theorem \ref{Pal-Vasu-P6-thm-003}}]
	Clearly, the upper bound of $K^p(B_{\ell^n_q},T,\lambda) $ can be obtained from Theorem \ref{Pal-Vasu-P6-thm-002}. To establish the lower bound we need an important result by Kwapie\'{n} \cite{kwapien-1968} which generalizes the well-known Grothendieck's theorem (for $s=2$) as well as Littlewood's theorem (for $T$ being the inclusion $\ell_1 \hookrightarrow \ell_{4/3}$): Every operator $T: \ell_1 \rightarrow \ell_s$ is $(r,1)$-summing, where 
	\begin{equation*}
		\frac{1}{r}=1-\bigg|\frac{1}{s}-\frac{1}{2}\bigg|.
	\end{equation*}
	Therefore, applying Theorem \ref{Pal-Vasu-P5-thm-B}(2) and by using the fact that $\ell_s$ is $\max\{2,s\}$-concave, we obtain our desired lower estimate. This completes the proof.
\end{pf}

\begin{pf}[{\bf Proof of Theorem \ref{Pal-Vasu-P6-thm-004}}]
	The lower bound of $A_p(B_{\ell^n_q},\ell_s\hookrightarrow\ell_t,\lambda)$ can be deduced from Lemma \ref{Pal-Vasu-P6-lem-001} and Theorem \ref{Pal-Vasu-P6-thm-002}. To establish the upper bound, we divide the proof into two different cases.\\[1mm]
	{\it The case $1\leq s\leq 2$.} This case can be settled by use of the result proved by Allu {\it et. al.} \cite[Theorem 2.2]{Subhadip-Vasu-P3}. In fact the authors \cite{Subhadip-Vasu-P3} have proved the following upper bound: For any $1\leq p<\infty$ and $\lambda \geq 1$ we have
	\begin{equation}\label{Pal-Vasu-P6-eqn-009}
		A_p(B_{\ell^n_q},X,\lambda) \prec 	\frac{(\log n)^{1-\frac{1}{\min\{q,2\}}}}{n^{\frac{1}{p}+\frac{1}{\max\{q,2\}}-\frac{1}{2}}}.
	\end{equation}
	It is worth noting that for any bounded operator $T$ we have 
	\begin{equation}\label{Pal-Vasu-P6-eqn-010}
		A_p(B_{\ell^n_q},T,\lambda)\leq A_p(B_{\ell^n_q},\mathbb{C},\lambda).
	\end{equation}
	In view of \eqref{Pal-Vasu-P6-eqn-009} and \eqref{Pal-Vasu-P6-eqn-010}, we obtain the desired upper bound for $A_p(B_{\ell^n_q},\ell_s\hookrightarrow\ell_t,\lambda)$.\\[1mm]
	{\it The case $2\leq s<\infty$.} Following the techniques from the Theorem \ref{Pal-Vasu-P6-thm-001}, we deal the upper bound in this case. Let $r=(r_1,\dots,r_n)\in \mathbb{R}^n_{\geq 0}$ be such that 
	\begin{equation*}
		\sum_{j=1}^{n}\norm{T(x_j(Q))}^p_{\ell^n_t}r_j \leq \sup_{z\in B_{\ell^n_{q}}}\norm{Q(z)}^p_{\ell^n_s}
	\end{equation*}
	for every $1$-homogeneous polynomial $Q \in \mathcal{P}(^1 B_{\ell^n_q}, B_{\ell^n_s})$, where $T:\ell^n_s \hookrightarrow \ell^n_t$ with $2\leq s\leq \infty$.
	Therefore, we have
	\begin{equation}\label{Pal-Vasu-P6-eqn-011}
		\sum_{j=1}^{n}r_j\frac{1}{n^{\frac{p}{q}}}\leq \sum_{j=1}^{n}r^p_j\frac{\norm{e_j}^p}{n^{\frac{p}{q}}}\leq \sup_{z\in B_{\ell^n_{q}}} \norm{\sum_{j=1}^{n}z_je_j}^p_{\ell^n_s} \leq \sup_{z\in B_{\ell^n_{q}}} \norm{z}^p_{\ell^n_s}.
	\end{equation}
	We observe that for $q\leq s$ we have $\sup_{z\in B_{\ell^n_{q}}}\norm{z}_{s}\leq 1$. Hence, in view of \eqref{Pal-Vasu-P6-eqn-011} we obtain
	\begin{equation*}
		\frac{1}{n}\sum_{j=1}^{n}r_j\leq \frac{1}{n^{1-\frac{p}{q}}}.
	\end{equation*}
	Therefore, it follows that 
	\begin{equation}\label{Pal-Vasu-P6-eqn-012}
		A^1_{p}(B_{\ell^n_q},\ell_s \hookrightarrow \ell_t) \leq \frac{1}{n^{1-\frac{p}{q}}}.
	\end{equation}
	On the other hand, if $q>s$ then we have 
	\begin{equation*}
		\sup_{z\in B_{\ell^n_{q}}}\norm{z}_{\ell^n_s}\leq n^{\frac{1}{s}-\frac{1}{q}}.
	\end{equation*}
	So, by the virtue of \eqref{Pal-Vasu-P6-eqn-011}, we obtain 
	\begin{equation*}
		\frac{1}{n}\sum_{j=1}^{n}r_j\leq \frac{1}{n^{1-\frac{p}{s}}}.
	\end{equation*}
	Consequently, it shows that 
	\begin{equation}\label{Pal-Vasu-P6-eqn-013}
		A^1_{p}(B_{\ell^n_q},\ell_s \hookrightarrow \ell_t) \leq \frac{1}{n^{1-\frac{p}{s}}}.
	\end{equation}
	From \eqref{Pal-Vasu-P6-eqn-012} and \eqref{Pal-Vasu-P6-eqn-013} and by using the following simple observation 
	\begin{equation*}
		A^1_{p}(B_{\ell^n_q},\ell_s \hookrightarrow \ell_t,\lambda)=\lambda^{1/p}	A^1_{p}(B_{\ell^n_q},\ell_s \hookrightarrow \ell_t),
	\end{equation*}
	we obtain
	$$
	A^1_{p}(B_{\ell^n_q},\ell_s \hookrightarrow \ell_t,\lambda) \leq  \left\{\begin{array}{ll}
		\frac{\lambda^{\frac{1}{p}}}{n^{1-\frac{p}{q}}}, & \mbox{for \, $q\leq s$}\\[5mm] 
		
		\frac{\lambda^{\frac{1}{p}}}{n^{1-\frac{p}{s}}}, & \mbox{for \, $q>s$}.
	\end{array}\right.
	$$
	Therefore, the final conclusion follows from the fact that
	\begin{equation*}
		A_{p}(B_{\ell^n_q},\ell_s \hookrightarrow  \ell_t,\lambda)\leq  A^1_{p}(B_{\ell^n_q},\ell_s \hookrightarrow \ell_t,\lambda).
	\end{equation*}
	This completes the proof.
\end{pf}

\begin{pf}[{\bf Proof of Theorem \ref{Pal-Vasu-P6-thm-005}}]
	We observe that the upper bound for $A_p(B_{\ell^n_q},T,\lambda)$ follows from Theorem \ref{Pal-Vasu-P6-thm-004}. On the other hand, the lower bound is just a consequence of Theorem \ref{Pal-Vasu-P6-thm-003} and Lemma \ref{Pal-Vasu-P6-lem-001}. This completes the proof.
\end{pf}
\vspace{3mm}

\noindent\textbf{Compliance of Ethical Standards:}\\

\noindent\textbf{Conflict of interest.} The authors declare that there is no conflict  of interest regarding the publication of this paper.
\vspace{1.5mm}

\noindent\textbf{Data availability statement.}  Data sharing is not applicable to this article as no datasets were generated or analyzed during the current study.
\vspace{1.5mm}

\noindent\textbf{Authors contributions.} Both the authors have made equal contributions in reading, writing, and preparing the manuscript.\\

\noindent\textbf{Acknowledgment:} 
The research of the second named author is supported by DST-INSPIRE Fellowship (IF 190721),  New Delhi, India.


\begin{thebibliography}{99}
	
		\bibitem{aizn-2000a} {\sc L. Aizenberg}, Multidimensional analogues of Bohr's theorem on power series, \textit{Proc. Amer. Math. Soc.} {\bf 128} (2000), 1147--1155.
	
	\bibitem{aizn-2000b} {\sc L. Aizenberg, A. Aytuna}, and {\sc P. Djakov}, An abstract approach to Bohr's phenomenon, {\it Proc. Amer. Math. Soc.} {\bf 128} (2000), 2611--2619.
	
	
	\bibitem{aizenberg-2001} {\sc L. Aizenberg, A. Aytuna}  and {\sc P. Djakov}, Generalization of theorem on Bohr for bases in spaces of holomorphic functions of several complex variables, 
	{\it J. Math. Anal. Appl.} {\bf  258} (2001), 429--447.
	
	
	
	
	
	
	
	
	\bibitem{Subhadip-Vasu-P3} {\sc V. Allu, H. Halder,} and {\sc S. Pal}, Multidimensional Bohr radii for vector-valued holomorphic functions, see https://arxiv.org/abs/2308.07825 (2023).
	
	\bibitem{Subhadip-Vasu-P4-forum} {\sc V. Allu, H. Halder,} and {\sc S. Pal}, Arithmetic Bohr radius for the Minkowski space, {\it Forum Math.} DOI: https://doi.org/10.1515/forum-2023-0425 (2024).
	
	\bibitem{Subhadip-Vasu-P5} {\sc V. Allu} and {\sc S. Pal}, On multidimensional Bohr radii for Banach spaces, see https://arxiv.org/abs/2406.19865 (2024).
	
	
	
	
	
	
	
	
	
	\bibitem{Bennett-1973} {\sc G. Bennett}, Inclusion mappings between $\ell^p$ spaces, {\it J. Funct. Anal.} {\bf 13} (1973), 20--27.
	
	\bibitem{Blasco-OTAA-2010} {\sc O. Blasco}, The Bohr radius of a Banach space, {\it In Vector measures, integration and related topics, 5964, Oper. Theory Adv. Appl., 201, Birkh{\"a}user Verlag, Basel, 2010.}
	
	\bibitem{Blasco-Collect-2017} {\sc O. Blasco}, The $p$-Bohr radius of a Banach space, {\it Collect. Math.} {\bf 68} (2017), 87--100.
	
	
	\bibitem{boas-1997} {\sc H. P. Boas} and {\sc D. Khavinson}, Bohr's power series theorem in several variables, {\it Proc. Amer. Math. Soc.}  {\bf 125} (1997), 2975--2979.
	
	\bibitem{boas-2000} {\sc H. P. Boas}, Majorant Series,  {\it J. Korean Math. Soc.}  {\bf 37} (2000), 321--337.
	
	\bibitem{Bohr-1914} {\sc H. Bohr}, A theorem concerning power series,  {\it Proc. Lond. Math. Soc.} s2-13 (1914), 1--5.
	
	
	\bibitem{bombieri-1962} {\sc E. Bombieri}, Sopra un teorema di H. Bohr e G. Ricci sulle funzioni maggioranti delle serie di potenze, {\it Boll. Un. Mat. Ital.} {\bf 17} (1962), 276--282.
	
	\bibitem{bombieri-2004} {\sc E. Bombieri} and {\sc J. Bourgain}, A remark on Bohr's inequality, {\it Internat. Math. Res. Notices} {\bf 80} (2004), 4307--4330.
	
	\bibitem{Carl-1974} {\sc B. Carl}, Absolut-$(p,1)$-summierende identische Operatoren von $\ell_u$ in $\ell_v$, {\it Math. Nachr.} {\bf 63} (1974), 353--360.
	
	
	
	\bibitem{defant-2006} {\sc A. Defant} and {\sc L. Frerick}, A logarithmic lower bound for multi-dimenional Bohr radii, {\it Israel J. Math.} {\bf 152} (2006), 17--28.
	
	
	
	\bibitem{defant-2003} {\sc A. Defant, D. Garc\'{i}a}, and {\sc  M. Maestre}, Bohr power series theorem and local Banach space theory, {\it J. Reine Angew. Math.} {\bf 557} (2003), 173--197.
	
		\bibitem{defant-2004} {\sc A. Defant, D. Garc\'{i}a}, and {\sc  M. Maestre}, Estimates for the first and second Bohr radii of Reinhardt domains, {\it J. Appr. Theory} {\bf 128} (2004), 53--68.
	
	\bibitem{defant-book} {\sc A. Defant, D. García,  M. Maestre}, and {\sc P. Sevilla-Peris}, \textit{Dirichlet Series and Holomorphic Functions in High Dimensions}, New Mathematical Monographs: {\bf 37}, Cambridge University Press, Cambridge, (2019).
	
	\bibitem{defant-QJM-2008} {\sc A. Defant,  M. Maestre}, and {\sc  C. Prengel}, The arithmetic Bohr radius, {\it Q. J. Math.} {\bf 59} (2008), 189--205.
	
	\bibitem{defant-Domain-J angew-2009} {\sc A. Defant,  M. Maestre}, and {\sc  C. Prengel}, Domains of convergence for monomial expansions of holomorphic functions in infinitely many variables, {\it J. Reine Angew. Math.} {\bf 634} (2009), 13--49.
	
	\bibitem{defant-adv-math-2012} {\sc A. Defant,  M. Maestre} and {\sc U. Schwarting}, Bohr radii of vector valued holomorphic functions, {\it Adv. Math.} {\bf 231} (2012), 2837--2857.
	
	\bibitem{diestel-Abs-sum-book} {\sc J. Diestel, H. Jarchow}, and {\sc A. Tonge}, {\it Absolutely Summing Operators}, in: Cambridge Studies in Advance Mathematics {\bf 43}, Cambridge University Press, Cambridge, (1995).
	
	
	\bibitem{Dixon & BLMS & 1995} {\sc P. G. Dixon}, Banach algebras satisfying the non-unital von Neumann inequality, \textit{Bull. Lond. Math. Soc.} \textbf{27} (4) (1995), 359--362.
	
	\bibitem{Dineen-Timoney-1989} {\sc S. Dineen} and {\sc R. M. Timoney}, Absolute bases, tensor products and a theorem of Bohr, \textit{Studia Math.} \textbf{94} (1989), 227--234.
	
	
	\bibitem{Djakov & Ramanujan & J. Anal & 2000} {\sc P. B. Djakov} and {\sc M. S. Ramanujan}, A remark on Bohr's theorem and its generalizations, \textit{J. Anal.} \textbf{8} (2000), 65--77.
	
	
	
	


	
	\bibitem{kumar-2023} {\sc S. Kumar}, On the multidimensional Bohr radius, {\it Proc. Amer. Math. Soc.} {\bf 151} (2023), 2001--2009.
	
	\bibitem{kumar-2023-arxiv} {\sc S. Kumar} and {\sc R. Manna}, Multi-dimensional Bohr radii of Banach space valued holomorphic functions, see https://arxiv.org/pdf/2303.17416 (2023).
	
	\bibitem{kwapien-1968} {\sc S. Kwapie\'{n}}, Some remarks on $(p,q)$-absolutely summing operators in $\ell_p$-spaces, {\it Stud. Math.} {\bf 29} (1968), 327--337.
	
	\bibitem{lindenstrauss-book-I} {\sc J. Lindenstrauss} and {\sc L. Tzafriri}, {\it Classical Banach Spaces I: Sequence Spaces}, in: Ergebnisse der Mathematik und ihrer Grenzgebiete, {\bf 92}, Springer-Verlag, Berlin-Heidelberg-New York, (1977).
	
	\bibitem{lindenstrauss-book-II} {\sc J. Lindenstrauss} and {\sc L. Tzafriri}, {\it Classical Banach Spaces II: Function Spaces}, in: Ergebnisse der Mathematik und ihrer Grenzgebiete, {\bf 97}, Springer-Verlag, Berlin-Heidelberg-New York, (1979).
	
	
	\bibitem{littlewood-1930} {\sc J. E. Littlewood}, On bounded bilinear forms in an infinite number of variables, {\it Q. J. Math.} {\bf 1} (1930), 164--174.
	
	\bibitem{paulsen-2002} {\sc V. I. Paulsen, G. Popescu} and {\sc D. Singh}, On Bohr's inequality, {\it Proc. Lond. Math. Soc.} s3-85 (2002), 493--512.
	
	
	
	
	\bibitem{sidon-1927} {\sc S. Sidon}, Uber einen satz von Hernn Bohr, {\it Math. Zeit.}  {\bf 26} (1927),  731--732.
	
	\bibitem{tomic-1962} {\sc M. Tomic}, Sur un theoreme de H. Bohr,  {\it Math. Scand.}  {\bf 11} (1962), 103--106.
\end{thebibliography}
\end{document}